\documentclass[a4paper,12pt]{amsart}

\usepackage{geometry}
\usepackage[comma,numbers]{natbib}
\usepackage[inline]{enumitem}
\usepackage{graphicx}
\usepackage[textfont=it,tableposition=above,font=small]{caption}		

\usepackage{framed}					
\usepackage[foot]{amsaddr}	
\usepackage{url}						

\usepackage{xargs}					
\usepackage{amssymb}				

\usepackage{silence}
\usepackage{todonotes}

\usepackage[ruled,vlined]{algorithm2e}

\usepackage{multirow}										
\usepackage{arydshln}										
\usepackage{psfrag}
\usepackage{mathtools}									

\usepackage{color}
\usepackage{xcolor}
\RequirePackage[colorlinks,citecolor=blue,urlcolor=blue]{hyperref}
\usepackage[norefs,nocites]{refcheck}

\usepackage{randtext}

\newtheorem{theorem}{Theorem}
\newtheorem{lemma}[theorem]{Lemma}

\theoremstyle{definition}

\newtheorem{example}{Example}

\newcommand{\twofigb}{.4}  

\DeclareMathOperator{\intr}{int}					
\DeclareMathOperator{\cl}{cl}   					
\DeclareMathOperator{\bd}{bd}                       
\DeclareMathOperator{\relint}{relint}				
\DeclareMathOperator{\relbd}{relbd} 				
\DeclareMathOperator{\relcl}{relcl}					
\DeclareMathOperator{\affOp}{aff}
\newcommand{\aff}[1]{\affOp\left(#1\right)}

\newcommand{\R}{\mathbb R}							
\newcommand{\Support}[1]{\mathrm{supp}\left(#1\right)}  
\newcommand{\Meas}[1][\R^d]{\mathcal M\left(#1\right)}
\newcommand{\I}[1]{\mathbb{I}\left[{#1}\right]}		

\newcommand{\Damu}[1][\mu]{D_\alpha(#1)}
\newcommand{\Uamu}[1][\mu]{U_\alpha(#1)}
\newcommand{\D}{D}			
\newcommand{\depth}[2]{\D\left(#1;#2\right)}
\newcommand{\half}{\mathcal{H}}
\newcommand{\flag}{\mathcal{F}}

\title[Partial Reconstruction of Measures from Halfspace Depth]{Partial Reconstruction of Measures from Halfspace Depth}

\author{Petra Laketa}
\author{Stanislav Nagy}
\email{\randomize{nagy@karlin.mff.cuni.cz}}

\address{
	Charles University,
	Faculty of Mathematics and Physics,
	Prague, Czech Republic
}
	
\date{\today}

\begin{document}
\begin{abstract}
The halfspace depth of a $d$-dimensional point $x$ with respect to a finite (or probability) Borel measure $\mu$ in $\R^d$ is defined as the infimum of the $\mu$-masses of all closed halfspaces containing $x$. A natural question is whether the halfspace depth, as a function of $x \in \R^d$, determines the measure $\mu$ completely. In general, it turns out that this is not the case, and it is possible for two different measures to have the same halfspace depth function everywhere in $\R^d$. In this paper we show that despite this negative result, one can still obtain a substantial amount of information on the support and the location of the mass of $\mu$ from its halfspace depth. We illustrate our partial reconstruction procedure in an example of a non-trivial bivariate probability distribution whose atomic part is determined successfully from its halfspace depth.
\end{abstract}

\maketitle


%
%
%

\section{The Depth Characterization/Reconstruction Problem}

Let $x$ be a point in the $d$-dimensional Euclidean space $\R^d$ and let $\mu$ be a finite Borel measure in $\R^d$. We write $\half$ for the collection of all closed halfspaces\footnote{A halfspace is one of the two regions determined by a hyperplane in $\R^d$; any halfspace can be written as a set $\left\{ y \in \R^d \colon \left\langle y, u \right\rangle \leq c \right\}$ for some $c \in \R$ and $u \in \R^d \setminus \left\{ 0 \right\}$.} in $\R^d$ and $\half(x)$ for the subset of those halfspaces from $\half$ that contain $x$ in their boundary hyperplane. The \emph{halfspace depth} (or \emph{Tukey depth}) of the point $x$ with respect to $\mu$ is defined as
    \begin{equation}    \label{halfspace depth}
    \D\left(x;\mu\right) = \inf_{H\in\half(x)} \mu(H).  
    \end{equation}
The history of the halfspace depth in statistics goes back to the 1970s \cite{Tukey1975}. The halfspace depth plays an important role in the theory and practice of nonparametric inference of multivariate data; for many references see \cite{Liu_etal1999,Nagy_etal2019,Zuo_Serfling2000}.

The depth~\eqref{halfspace depth} was originally designed to serve as a multivariate generalization of the quantile function. As such, it is desirable that just as the quantile function in $\R$, the depth function $x \mapsto \D(x;\mu)$ in $\R^d$ characterizes the underlying measure $\mu$ uniquely, and $\mu$ is straightforward to be retrieved from its depth. The question whether the last two properties are valid for $\D$ are known as the \emph{halfspace depth characterization and reconstruction problems}. They both turned out not to have an easy answer. In fact, only the recent progress in the theory of the halfspace depth gave the first definite solutions to some of these problems.
    
In \cite{Nagy2019b}, the general characterization question for the halfspace depth was answered in the negative, by giving examples of different probability distributions in $\R^d$ with $d \geq 2$ with identical halfspace depth functions. On the other hand, several authors have obtained also partial positive answers to the characterization problem; for a recent overview of that work see \cite{Nagy2020c}. Only three types of distributions are known to be completely characterized by their halfspace depth functions: \begin{enumerate*}[label=(\roman*)] \item univariate measures, in which case the depth~\eqref{halfspace depth} is just a simple transform of the distribution function of $\mu$; \item atomic measures with finitely many atoms (which we subsequently call \emph{finitely atomic measures} for brevity) in $\R^d$ \cite{Struyf_Rousseeuw1999,Laketa_Nagy2021}; and \item measures that possess all Dupin floating bodies\footnote{A Borel measure $\mu$ on $\R^d$ is said to possess all Dupin floating bodies if each tangent halfspace to the halfspace depth upper level set $\left\{ x \in \R^d \colon \D(x;\mu) \geq \alpha \right\}$ is of $\mu$-mass exactly $\alpha$, for all $\alpha \geq 0$.} \cite{Nagy_etal2019}.\end{enumerate*}

In this contribution we revisit the halfspace depth reconstruction problem. We pursue a general approach, and do not restrict only to atomic measures, or to measures with densities. Our results are valid for any finite (or probability) Borel measure $\mu$ in $\R^d$. As the first step in addressing the reconstruction problem, our intention is to identify the support and the location of the atoms of $\mu$, based on its depth. We will see at the end of this note that without additional assumptions, neither of these problems is possible to be resolved. We, however, prove several positive results. 

We begin by introducing the necessary mathematical background in Section~\ref{section:preliminaries}. In Section~\ref{section:main} we state our main theorem; a detailed proof of that theorem is given in the Appendix. We show that \begin{enumerate*}[label=(\roman*)] \item the support of the measure $\mu$ must be concentrated only in the boundaries of the level sets of its halfspace depth; \item each atom of $\mu$ is an extreme point of the corresponding (closed and convex) upper level sets of the halfspace depth; and \item each atom of $\mu$ induces a jump in the halfspace depth function on the line passing through that atom. \end{enumerate*} These advances enable us to identify the location of the atoms of $\mu$, at least in simpler scenarios. We illustrate this in Section~\ref{section:examples}, where we give an example of a non-trivial bivariate probability measure $\mu$ whose atomic part we are able to determine from its depth. We conclude by giving an example of two measures whose depth functions are the same, yet both their supports and the location of their atoms differ.

\section{Preliminaries: Flag Halfspaces and Central Regions}    \label{section:preliminaries}

\subsubsection*{Notations.}

When writing simply a subspace of $\R^d$ we always mean an affine subspace, that is the set $a + L = \left\{ a + x \in \R^d \colon x \in L \right\}$ for $a \in \R^d$ and $L$ a linear subspace of $\R^d$. The intersection of all subspaces in $\R^d$ that contain a set $A \subseteq \R^d$ is called the affine hull of $A$, and denoted by $\aff{A}$. It is the smallest subspace that contains $A$. The affine hull $\aff{\{x,y\}}$ of two different points $x, y \in \R^d$ is the infinite line passing through both $x$ and $y$; another example of a subspace is any hyperplane in $\R^d$.  

For a set $A\subseteq \R^d$ we write $\intr(A)$, $\cl(A)$ and $\bd(A)$ to denote the interior, closure, and boundary of $A$, respectively. The interior, closure, and boundary of a set $B\subseteq A$ when considered only as a subset of a subspace $A \subseteq \R^d$ are denoted by $\intr_A(B)$, $\cl_A(B)$ and $\bd_A(B)$, respectively. For two different points $x,y\in\R^d$, $x \ne y$, we denote by $L(x,y)$ the interior of the line segment between $x$ and $y$ when considered inside the infinite line $\aff{\{x,y\}}$. In other words, $L(x,y)$ is the open line segment between $x$ and $y$. In the special case of $A=\aff{B}$ we write $\relint(B)=\intr_A(B)$, $\relbd(B)=\bd_A(B)$ and $\relcl(B)=\cl_A(B)$ to denote the relative interior, relative boundary, and relative closure of $B$, respectively. For instance, $\relbd(L(x,y)) = \{x, y \}$ and $L(x,y) = \relint(L(x,y))$, but $\intr(L(x,y)) = \emptyset$ if $d>1$.

We write $\Meas$ for the collection of all finite Borel measures in $\R^d$. For a subspace $A\subseteq \R^d$ and $\mu\in\Meas$ we write $\mu|_A$ to denote the measure obtained by restricting $\mu$ to the subspace $A$, that is the finite Borel measure given by $\mu|_A(B) = \mu(B \cap A)$ for any Borel set $B \subseteq \R^d$. By $\Support{\mu}$ we mean the support of $\mu\in\Meas$, which is the smallest closed subset of $\R^d$ of full $\mu$-mass.

\subsection{Minimizing halfspaces and flag halfspaces}

For $\mu\in\Meas$ and $x \in \R^d$ we call $H \in \half(x)$ a \emph{minimizing halfspace} of $\mu$ at $x$ if $\mu(H) = \D\left(x;\mu\right)$. For $d = 1$ a minimizing halfspace always trivially exists. It also exists if $\mu$ is smooth in the sense that $\mu(\bd(H)) = 0$ for all $H \in \half(x)$, or if $\mu \in \Meas$ is finitely atomic. In general, however, the infimum in~\eqref{halfspace depth} does not have to be attained. We give a simple example.

\begin{figure}[htpb]
\begin{center}
\includegraphics[width=\twofigb\textwidth]{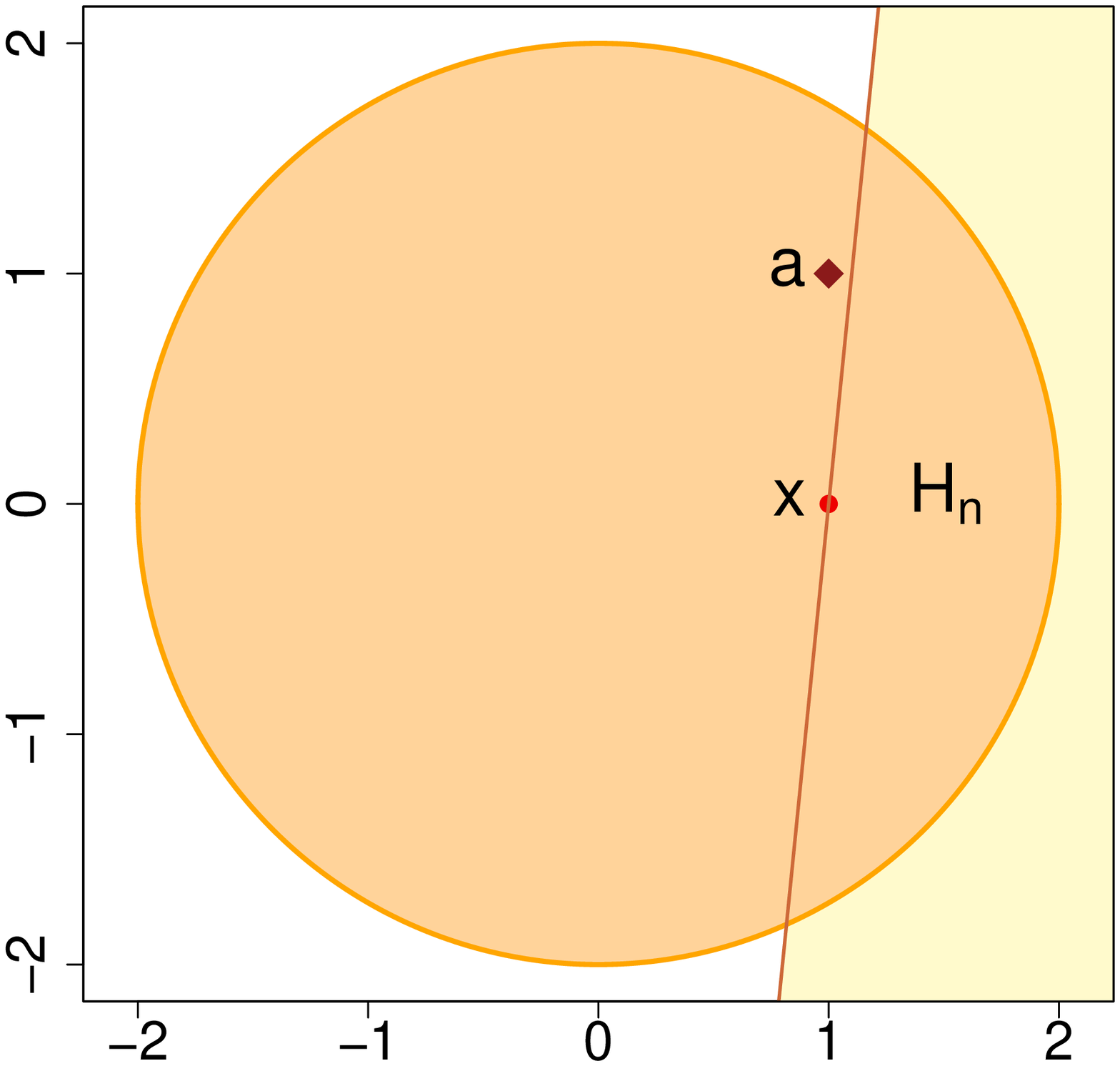} \qquad \includegraphics[width=\twofigb\textwidth]{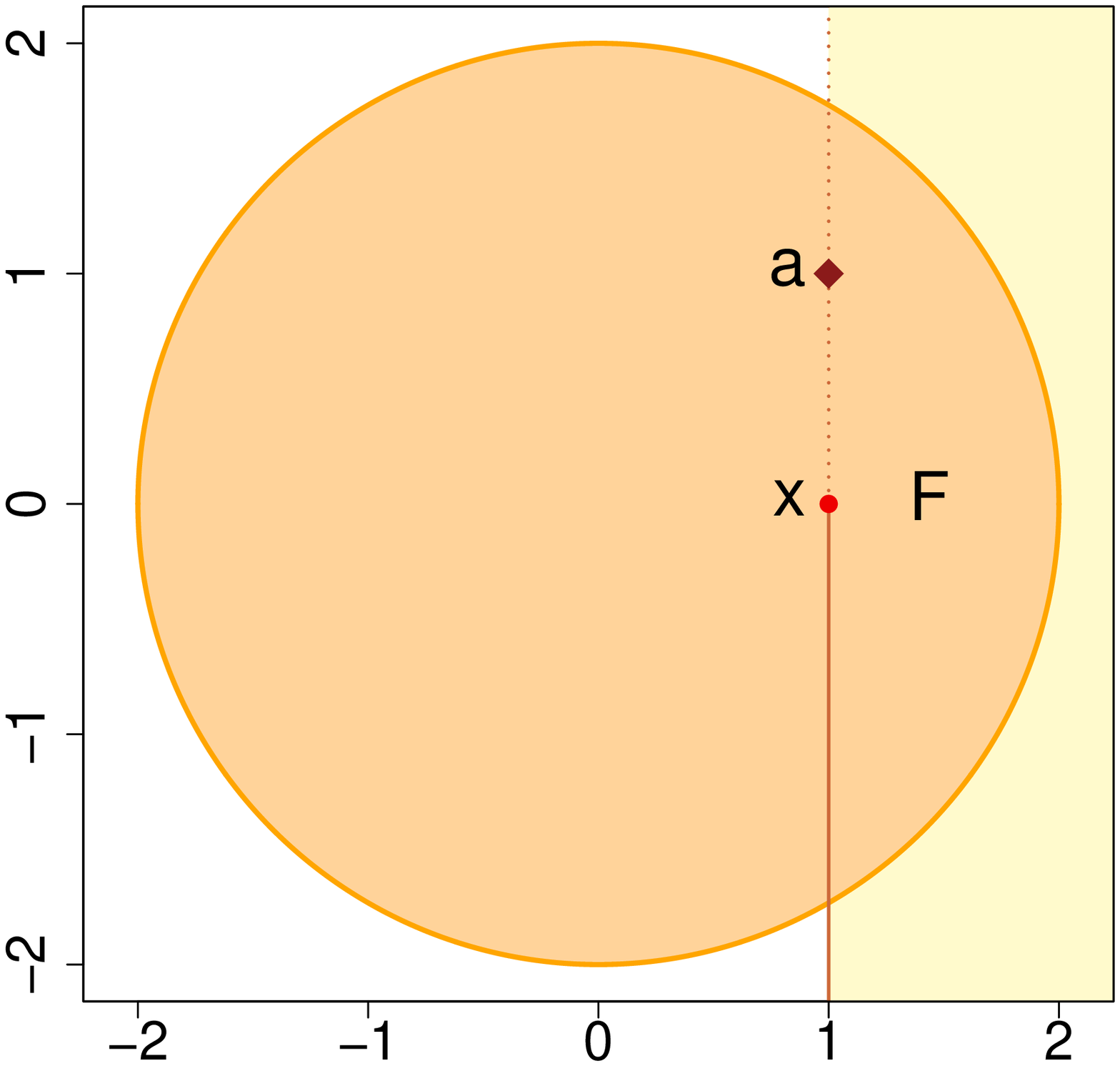}%
\end{center}
\caption{The support of $\mu \in \Meas[\R^2]$ from Example~\ref{example:flag} (colored disk) and its unique atom $a$ (diamond). No minimizing halfspace of $\mu$ at $x = (1,0)\in\R^2$ exists. On the left hand panel we see a halfspace $H_n \in \half(x)$ whose $\mu$-mass is almost $\D(x;\mu)$. The halfspace $H_n$ does not contain $a$. On the right hand panel the unique minimizing flag halfspace $F \in \flag(x)$ of $\mu$ at $x$ is displayed.}
\label{figure:flag halfspace}
\end{figure}

\begin{example}	\label{example:flag}
Take $\mu \in \Meas[\R^2]$ the sum of a uniform distribution on the disk $B = \left\{ x \in \R^2 \colon \left\Vert x \right\Vert \leq 2 \right\}$ and a Dirac measure at $a = (1,1) \in \R^2$. For $x = (1,0) \in \R^2$ no minimizing halfspace of $\mu$ at $x$ exists. As can be seen in Fig.~\ref{figure:flag halfspace}, the depth $\D(x;\mu)$ is approached by $\mu(H_{n})$ for a sequence of halfspaces $H_n\in\half(x)$ with inner normals $v_n = \left( \cos(-1/n), \sin(-1/n) \right)$ that converge $H_{v} \in \half(x)$ with inner normal $v = (1,0)$, yet $\D(x;\mu) = \lim_{n \to \infty} \mu(H_n) < \mu(H_{v})$.
\end{example}

For certain theoretical properties of the halfspace depth of $\mu$ to be valid, the existence of minimizing halfspaces appears to be crucial. As a way to alleviate the issue of their possible non-existence, in \cite{Pokorny_etal2022} a novel concept of the so-called flag halfspaces was introduced. A \emph{flag halfspace} $F$ centered at a point $x\in\R^d$ is defined as any set of the form 
    \begin{equation}\label{flag halfspace}
     F=\{x\} \cup \left( \bigcup_{i=1}^d \relint(H_i) \right).
    \end{equation}
In this formula, $H_d\in\half(x)$ and for each $i\in \{1,\dots,d-1\}$, $H_i$ stands for an $i$-dimensional halfspace inside the subspace $\relbd(H_{i+1})$ such that $x\in\relbd(H_i)$. The collection of all flag halfspaces in $\R^d$ centered at $x\in\R^d$ is denoted by $\flag(x)$. An example of a flag halfspace in $\R^2$ is displayed in the right hand panel of Fig.~\ref{figure:flag halfspace}. That flag halfspace is a union of an open halfplane $H_2$ (light-colored halfplane) whose boundary passes through $x$, a halfline (thick halfline) in the boundary line $\bd(H_2)$ starting at $x$, and the point $x$ itself.

The results derived the present paper lean on the following crucial observation, whose complete proof can be found in \cite[Theorem~2]{Pokorny_etal2022}.

\begin{lemma}\label{theorem:Pokorny}
For any $x\in\R^d$ and $\mu\in\Meas$ it holds true that
    \[  \D\left(x;\mu\right)=\min_{F\in\flag(x)}\mu(F). \]
In particular, there always exists $F\in\flag(x)$ such that $\mu(F)=\D\left(x;\mu\right)$.
\end{lemma}

Any flag halfspace $F \in \flag(x)$ from Lemma~\ref{theorem:Pokorny} that satisfies $\mu(F)=\D\left(x;\mu\right)$ is called a \emph{minimizing flag halfspace} of $\mu$ at $x$. This is because it minimizes the $\mu$-mass among all the flag halfspaces from $\flag(x)$. Lemma~\ref{theorem:Pokorny} tells us two important messages. First, the halfspace depth $\D(x;\mu)$ can be introduced also in terms of flag halfspaces instead of the usual closed halfspaces in~\eqref{halfspace depth}, and the two formulations are equivalent to each other. Second, in contrast to the usual minimizing halfspaces that do not exist at certain points $x \in \R^d$, according to Lemma~\ref{theorem:Pokorny} there always exists a minimizing flag halfspace of any $\mu$ at any $x$.

\subsection{Halfspace depth central regions}

The upper level sets of the halfspace depth function $\D(\cdot;\mu)$, given by  
    \begin{equation}	\label{central region}
	\Damu = \left\{ x \in \R^d \colon \depth{x}{\mu} \geq \alpha \right\} \mbox{ for }\alpha\geq 0,
	\end{equation}
play the important role of multivariate quantiles in depth statistics. The set $\Damu$ is called the \emph{central region} of $\mu$ at level $\alpha$. All central regions are known to be convex and closed. The sets~\eqref{central region} are clearly also nested, in the sense that $\Damu \subseteq \D_\beta(\mu)$ for $\beta\leq \alpha$. Besides~\eqref{central region}, another collection of depth-generated sets of interest considered in \cite{Laketa_Nagy2021b,Pokorny_etal2022} is 
    \begin{equation*}	
	\Uamu = \left\{ x \in \R^d \colon \depth{x}{\mu} > \alpha \right\} \mbox{ for }\alpha\geq 0.
	\end{equation*}
We conclude this collection of preliminaries with another result from \cite{Pokorny_etal2022}, which tells us that no set difference of the level sets $\Damu\setminus\Uamu$ can contain a relatively open subset of positive $\mu$-mass. That result lends an insight into the properties of the support of $\mu$, based on its depth function $\D\left(\cdot;\mu\right)$. It will be of great importance in the proof of our main result in Section~\ref{section:main}. The complete proof of the next technical lemma can be found in~\cite[Lemma~6]{Pokorny_etal2022}.

\begin{lemma}\label{main lemma for support}
Let $\mu \in \Meas$ and let $K \subset \R^d$ be a relatively open set of points of equal depth of $\mu$ that contains at least two points. Then $\mu(K) = 0$. 
\end{lemma}

\section{Main Result}   \label{section:main}

The preliminary Lemma~\ref{main lemma for support} hints that the mass of $\mu$ cannot be located in the interior of regions of constant depth. We refine and formalize that claim in the following Theorem~\ref{theorem:support}, which is the main result of the present work. 

In part~\ref{support2} of Theorem~\ref{theorem:support} we bound the support of $\mu\in\Meas$, based on the information available in its depth function $\D\left(\cdot;\mu\right)$. We do so by showing that $\mu$ may be supported only in the closure of the boundaries of the central regions $\Damu$. That is a generalization of a similar result, known to be valid in the special case of finitely atomic measures $\mu\in\Meas$ \cite{Laketa_Nagy2021,Liu_etal2020,Struyf_Rousseeuw1999}. In the latter situation, all central regions $\Damu$ are convex polytopes, there is only a finite number of different polytopes in the collection $\left\{ \Damu \colon \alpha \geq 0 \right\}$, and the atoms of $\mu$ must be located in the vertices of the polytopes from that collection. Nevertheless, not all vertices of $\Damu$ are atoms of $\mu$; an algorithmic procedure for the reconstruction of the atoms, and the determination of their $\mu$-masses, is given in \cite{Laketa_Nagy2021}. 

Extending the last observation about the possible location of atoms from finitely atomic measures to the general scenario, in part~\ref{support1} of Theorem~\ref{theorem:support} we show that all atoms of $\mu$ are contained in the extreme points\footnote{For a convex set $C \subset \R^d$, a face of $C$ is a convex subset $F \subseteq C$ such that $x, y \in C$ and $(x+y)/2 \in F$ implies $x, y \in F$. If $\{z\}$ is a face of $C$, then $z$ is called an extreme point of $C$.} of the central regions $\Damu$. Note that this indeed corresponds to the known theory for finitely atomic measures --- the extreme points of polytopes are exactly their vertices.

Our last observation in part~\ref{jump} of Theorem~\ref{theorem:support} is that each atom $x\in\R^d$ of $\mu$ induces a jump discontinuity in the halfspace depth, when considered on the straight line connecting any point of higher depth with $x$. This will be useful in detecting possible locations of atoms for general measures.
	
\begin{theorem}	\label{theorem:support}
Let $\mu\in\Meas$.
\begin{enumerate}[label=(\roman*)]
\item \label{support2} Let $A$ be a subspace of $\R^d$ that contains at least two points. Then 
	\[	\Support{\mu|_A} \subseteq \cl_A\left(\bigcup_{\alpha \geq 0}\bd_A\left(\Damu\cap A\right)\right),	\]
	In particular, for $A = \R^d$ we have
	\[	\Support{\mu} \subseteq \cl\left(\bigcup_{\alpha \geq 0}\bd\left(\Damu\right)\right).	\]
\item \label{support1} Each atom $x$ of $\mu$ with $\depth{x}{\mu} = \alpha$ is an extreme point of $\D_\beta(\mu)$ for all $\beta \in (\alpha - \mu(\{x\}), \alpha]$.
\item \label{jump} For any $x \in \R^d$ with $\depth{x}{\mu} = \alpha$, any $z\in\Uamu$, and any $y \in \R^d$ such that $x$ belongs to the open line segment $L(y,z)$ between $y$ and $z$, it holds true that
    \[	\D\left(y;\mu\right)\leq\D\left(x;\mu\right)-\mu(\{x\}).\]
\end{enumerate}
\end{theorem}

The proof of Theorem~\ref{theorem:support} is given in the Appendix. Theorem~\ref{theorem:support} sheds light on the support and the location of the atoms of a measure. Its part~\ref{support2} tells us that if a depth function $\D\left(\cdot;\mu\right)$ attains only at most countably many different values, and each level set $\Damu$ is a polytope, the mass of $\mu$ must be concentrated in the closure of the set of vertices of the level sets $\Damu$. A special case is, of course, the setup of finitely atomic measures treated in \cite{Struyf_Rousseeuw1999,Laketa_Nagy2021}.

\section{Examples}  \label{section:examples}

We conclude this note by giving two examples. Parts~\ref{support1} and~\ref{jump} of Theorem~\ref{theorem:support} show a way, at least in special situations, to locate the atomic parts of measures. We start by reconsidering our motivating Example~\ref{example:flag}. The distribution $\mu\in\Meas[\R^2]$ is not purely atomic, and can be shown not to possess Dupin floating bodies. Thus, it is currently unknown whether its depth function $\D\left(\cdot;\mu\right)$ determines $\mu$ uniquely. In our first example of this section we show how Theorem~\ref{theorem:support} recovers the position the atomic part of $\mu$. Then, in Example~\ref{example:Nagy} we argue that the general problem of determining the support, or the location of the atoms of $\mu \in \Meas$ from its halfspace depth is impossible to be solved without further restrictions.

\addtocounter{example}{-1}
\begin{example}[continued]
Suppose that in Example~\ref{example:flag} we have $\mu(\{a\}) = \delta$ for $\delta \in (0,1/2)$ small enough, and that the non-atomic part of $\mu$ is $\nu\in\Meas[\R^2]$ uniform on the disk $B$, with $\nu(B) = 1$. Hence, $\mu(\R^2) = \nu(B) + \mu(\{a\}) = 1 + \delta$. We first compute the halfspace depth function $\D\left(\cdot;\mu\right)$ of $\mu$, and then show how to use Theorem~\ref{theorem:support} to find the atom $a$ of $\mu$ from its depth. The computation of the depth function is done by means of determining all the central regions~\eqref{central region} at levels $\beta \geq 0$ of $\mu$. We denote $\alpha = \D(x;\mu)$, and split our argumentation into three situations according to the behavior of the regions $\D_{\beta}(\mu)$: \begin{enumerate*}[label=(\roman*)] \item $\beta < \alpha$ where $x$ is contained in the interior of $\D_{\beta}(\mu)$; \item $\beta \in (\alpha, \alpha + \delta]$ where $x$ lies in the boundary of $\D_{\beta}(\mu)$; and \item $\beta > \alpha + \delta$ where $\D_{\beta}(\mu)$ does not contain $x$.\end{enumerate*} First note that because $\nu$ is uniform on a unit disk, all non-empty depth regions $D_{\beta}(\nu)$ of $\nu$ are circular disks centered at the origin, and all the touching halfspaces\footnote{We say that $H\in\half$ is \emph{touching} $A\subset\R^d$ if $H\cap A\neq \emptyset$ and $\intr(H)\cap A=\emptyset$.} of $\D_{\beta}(\nu)$ carry $\nu$-mass exactly $\beta$.

\begin{figure}[htpb]
\begin{center}
\includegraphics[width=\twofigb\linewidth]{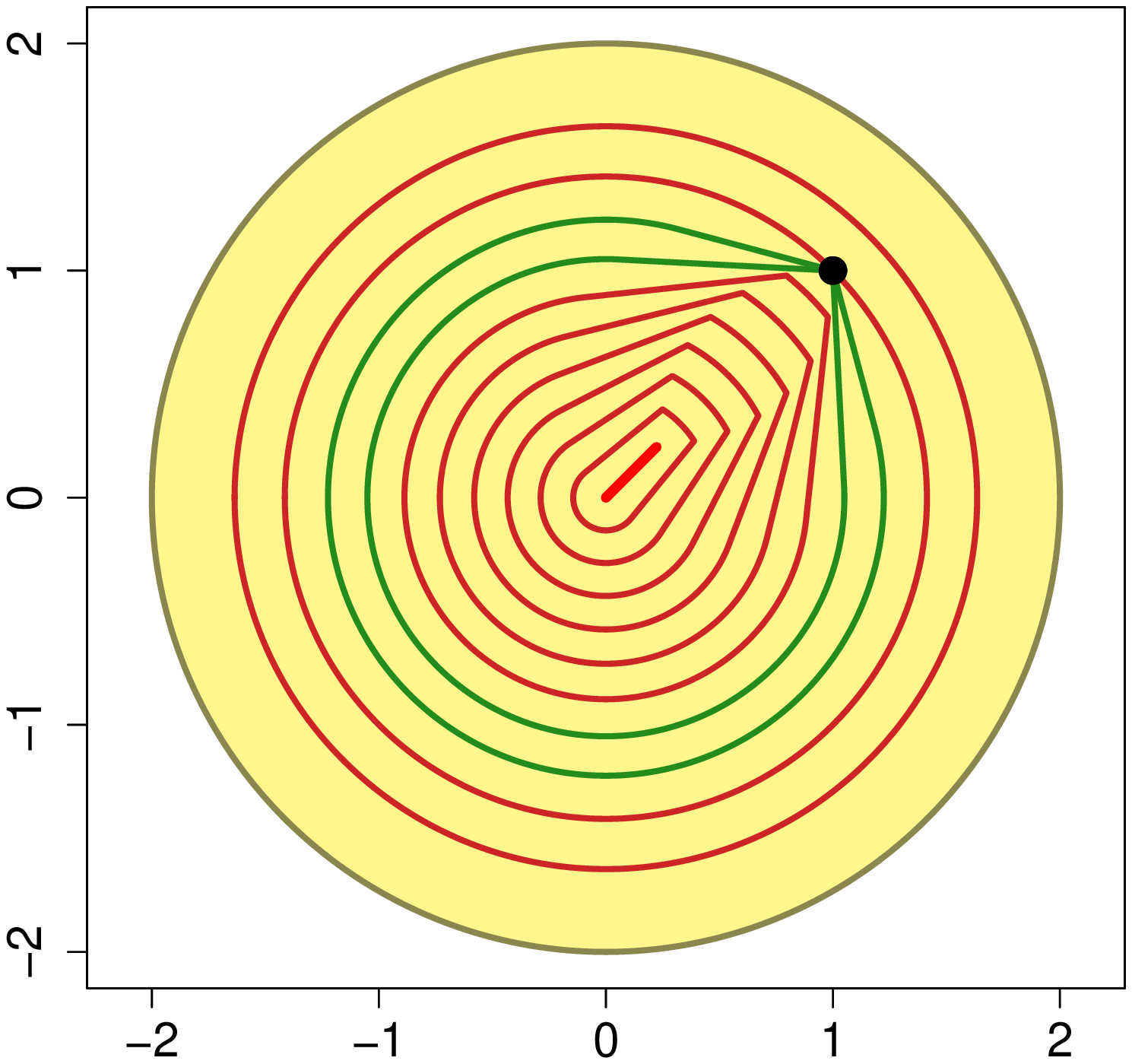} \qquad
\includegraphics[width=\twofigb\linewidth]{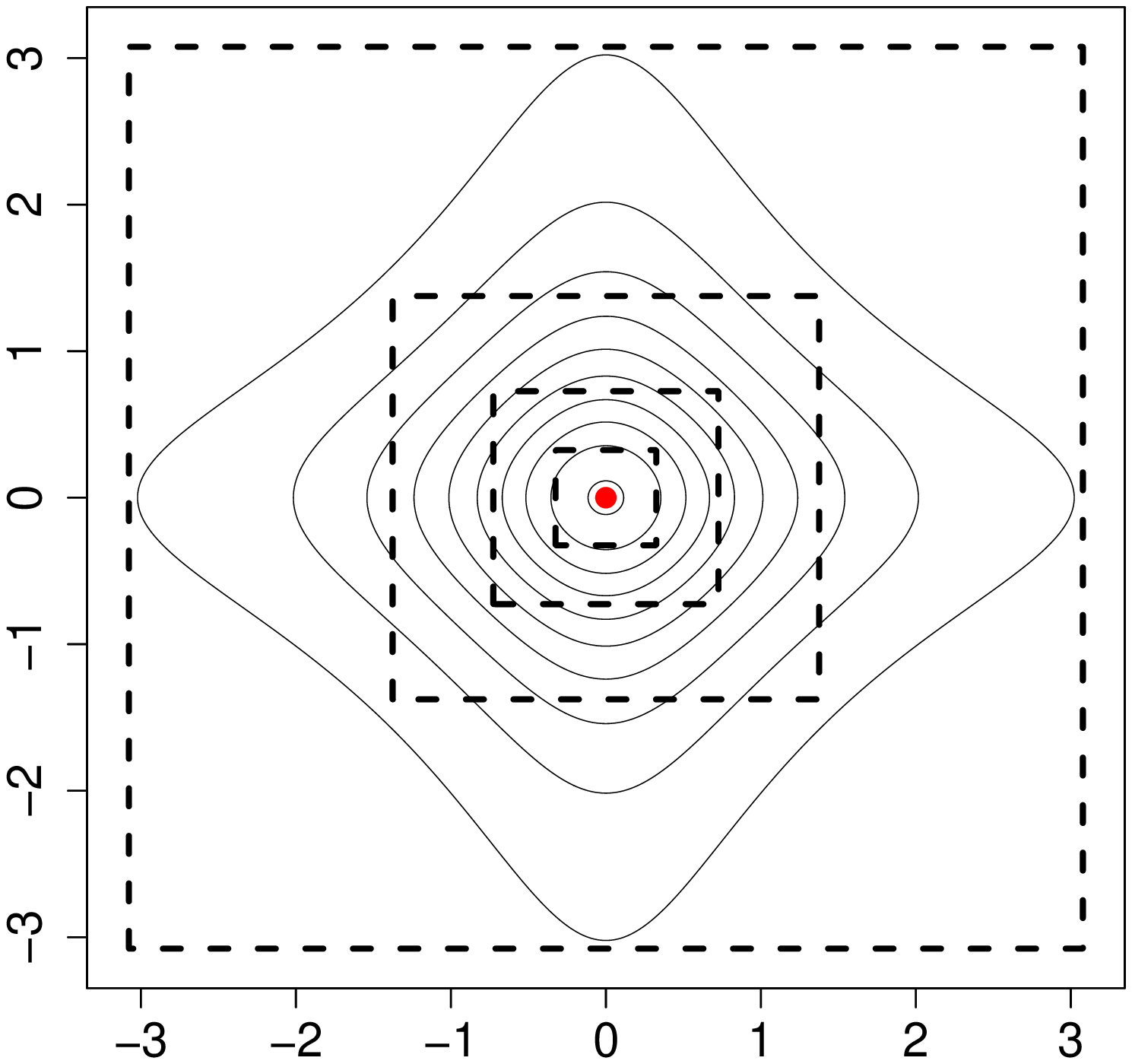}
\end{center}
\caption{Left panel: The measure $\mu$ from Example~\ref{example:flag} being the sum of a uniform distribution on a disk and a single atom at $a \in \R^2$ (black point) with $\delta = 1/10$, with several contours of its depth $\D\left(\cdot;\mu\right)$ (thick colored lines). The halfspace median is the red line segment in the middle of the plot. From the depth only, Theorem~\ref{theorem:support} allows us to determine the mass and the location of the atom. Two depth contours that share $a \in \R^2$ as an extreme point are plotted with boundaries in green. Right panel: Example~\ref{example:Nagy} with $d=2$. Several density contours of the measure $\mu \in \Meas[\R^2]$ (solid lines) and its atom (point at the origin), together with multiple contours of the corresponding depth $\D(\cdot;\mu) \equiv \D(\cdot; \nu)$ (dashed lines).}
\label{figure:atom}
\end{figure}

\smallskip\noindent
\textbf{Case I:} $\beta\leq \alpha$. For $\alpha = \depth{a}{\mu} = \depth{a}{\nu}$ we have that $\D_{\alpha}(\nu)$ is a disk centered at the origin containing $a$ on its boundary. Note that the halfspace depths of $\mu$ and $\nu$ remain the same outside $\D_{\alpha}(\nu)$, since the added atom $a$ does not lie in any minimizing halfspace of $x \notin \D_\alpha(\nu)$, so we have $\D_{\beta}(\mu)=\D_{\beta}(\nu)$ for all $\beta \leq \alpha$.

\smallskip\noindent
\textbf{Case II:} $\beta\in (\alpha,\alpha+\delta]$. We have $\depth{a}{\mu}=\alpha+\delta\geq \beta$, meaning that $a\in D_{\beta}(\mu)$. Because $\mu$ is obtained by adding mass to $\nu$, it must be $D_{\beta}(\nu)\subseteq D_{\beta}(\mu)$ and due to the convexity of the central regions \eqref{central region}, the convex hull $C$ of $D_{\beta}(\nu) \cup \{a\}$ must be contained in $D_{\beta}(\mu)$. Denote by $H \in \half(a)$ a touching halfspace of $D_{\beta}(\nu)$ that contains $a$ on its boundary. Then $\nu(H)=\beta$, and hence $\intr(H)\cap D_{\beta}(\mu)=\emptyset$. We obtain that $D_{\beta}(\mu)$ is equal to the convex hull of $D_{\beta}(\nu)\cup\{a\}$.

\smallskip\noindent
\textbf{Case III:} $\beta>\alpha+\delta$. In a manner similar to Case II one concludes that $D_{\beta}(\mu)$ is the convex hull of a circular disk $D_{\beta}(\nu)$ and $a$, intersected with the disk $D_{\beta-\delta}(\nu)$.

\smallskip

In order to complete the reconstruction of the atomic part of measure $\mu$ from Example~\ref{example:flag} based on its depth function, we present Lemma~\ref{lemma:Laketa}, which is a special case of a more general result (called the \emph{generalized inverse ray basis theorem}) whose complete proof can be found in \cite[Lemma~4]{Laketa_Nagy2021b}. 

\begin{lemma}	\label{lemma:Laketa}
Suppose that $\mu\in\Meas$, $\alpha > 0$, a point $x\notin \Damu$ and a face $F$ of $\Damu$ are given so that the relatively open line segment $L(x,y)$ does not intersect $\Damu$ for any $y \in F$. Then there exists a touching halfspace $H \in \half$ of $\Damu$ such that $\mu(\intr(H))\leq\alpha$, $x\in H$, and $F\subset\bd(H)$.
\end{lemma}

\smallskip\noindent
\emph{Reconstruction.} We now know the complete depth function $\D\left(\cdot;\mu\right)$ of $\mu$, see also Fig.~\ref{figure:atom}. From this depth only, we will locate the atoms of $\mu$ and their mass. The only point in $\R^2$ that is an extreme point of more than one depth region is certainly $a$, so that $a$ is the only possible candidate for an atom of $\mu$ by part~\ref{support1} of Theorem~\ref{theorem:support}. Take any $\beta\in(\alpha,\alpha+\delta)$. Then $D_{\beta}(\mu)$ is the convex hull of a circular disk and the point $a$ outside that disk, so its boundary contains a line segment $L(a,y_\beta)$ for $y_\beta\in \bd(D_{\beta}(\nu))$. Due to Lemma~\ref{lemma:Laketa}, there is a halfspace $H_{\beta} \in \half$ such that $L(a,y_\beta)\subset\bd(H_{\beta})$ and $\mu(\intr(H_{\beta}))\leq \beta < \alpha+\delta=\depth{a}{\mu} \leq \mu(H_\beta)$, the last inequality because $a \in H_\beta$. We obtain $\mu(\bd(H_{\beta}))\geq \alpha + \delta - \beta$. This is true for any $\beta \in (\alpha, \alpha + \delta)$, and for different $\beta_1, \beta_2 \in (\alpha,\alpha+\delta)$ we have $H_{\beta_1}\neq H_{\beta_2}$ with $x \in \bd(H_{\beta_i})$ and $\mu\left(\bd(H_{\beta_i})\right) \geq \alpha + \delta - \beta_i$, $i=1,2$. In conclusion, we obtain uncountably many different lines $\bd(H_\beta)$ of positive $\mu$-mass, all passing through $a$. That is possible only if $a$ is an atom of $\mu$, and $\mu(\{a\}) \geq \delta$. Theorem~\ref{theorem:support} again guarantees that $\mu(\{a\}) = \delta$ and that there is no other atom of $\mu$.
\end{example}

The complete Example~\ref{example:flag} gives a partial positive result toward the halfspace depth characterization problem, and promises methods allowing one to determine features of $\mu$ from its depth $\D\left(\cdot;\mu\right)$, at least for special sets of measures. The complete determination of the support or the atoms of $\mu$ from its depth is, however, a problem considerably more difficult, impossible to be solved in full generality. Follows an example of mutually singular measures\footnote{Recall that $\mu, \nu \in \Meas$ are called \emph{singular} if there is a Borel set $A \subset \R^d$ such that $\mu(A) = \nu(\R^d\setminus A) = 0$.} $\mu, \nu \in \Meas$ sharing the same depth function from \cite[Section~2.2]{Nagy2020c}. 

\begin{example}	\label{example:Nagy}
For $\mu_1 \in \Meas[\R^d]$ with independent Cauchy marginals and $\mu_2 \in \Meas[\R^d]$ the Dirac measure at the origin, define $\mu \in \Meas[\R^d]$ by the sum of $\mu_1$ and $\mu_2$ with weights $1/d$ and $1/2 - 1/(2d)$, respectively. The total mass of $\mu$ is hence $\mu\left(\R^d\right) = 1/2 + 1/(2d)$, and its support is $\R^d$. For the other distribution take $\nu \in \Meas[\R^d]$ the probability measure supported in the coordinate axes $A_i = \left\{ x = \left(x_1, \dots, x_d \right) \in \R^d \colon x_j = 0 \mbox{ for all }j \ne i \right\}$, $i=1,\dots,d$. The density $g$ of $\nu$ with respect to the one-dimensional Hausdorff measure on its support $\Support{\nu} = \bigcup_{i=1}^d A_i$ is given as a weighted sum of densities of univariate Cauchy distributions in $A_i$
	\[	g(x) = \frac{1}{d} \sum_{i=1}^d \frac{\I{x \in A_i}}{\pi(1+x_i^2)}	\quad\mbox{for }x = \left(x_1, \dots, x_d\right)\in\R^d.	\]
It can be shown \cite[Section~2.2]{Nagy2020c} that the depths of $\mu$ and $\nu$ coincide at all points $x = \left(x_1, \dots, x_d\right)$ in $\R^d$
	\[ 
	\depth{x}{\mu} = \depth{x}{\nu} = 
	\begin{cases}
	\frac{1}{d} \left( \frac{1}{2} - \frac{{\arctan(\max_{i=1,\dots,d} \left\vert x_i \right\vert)}}{\pi} \right) & \mbox{if }x \in \R^d \setminus \{ 0 \}, \\
	1/2 & \mbox{for }x = 0 \in \R^d.
	\end{cases}	
	\]
The two measures $\mu$ and $\nu$ are, however, singular as for $A = \Support{\nu}\setminus\{0\}$ we have $\mu(A) = \nu(\R^d \setminus A) = 0$. For an arbitrary finite Borel measure, it is therefore impossible to retrieve the full information about its support only from its depth function. For a visualization of measure $\mu$ and its halfspace depth see Fig.~\ref{figure:atom}.

The same example demonstrates that in general, also the location of the atoms of $\mu\in\Meas$, or even the number of them, cannot be recovered from the depth function $\D\left(\cdot;\mu\right)$ only --- the measure $\nu$ in Example~\ref{example:Nagy} does not contain any atoms, but $\mu$ has a single atom at its unique halfspace median (the smallest non-empty central region~\eqref{central region}). Because of the very special position of the atom of $\mu$, it is impossible to use our results from parts~\ref{support1} and~\ref{jump} of Theorem~\ref{theorem:support} to decide whether the origin is an atom of $\mu$, or not.
\end{example}

\section{Conclusion}

The halfspace depth has many applications, for example in classification or in nonparametric goodness-of-fit testing. However, in order to apply it properly, one needs to make sure that the measure $\mu$ is characterized by its halfspace depth function, so that we can use the halfspace depth to distinguish $\mu$ from other measures. For that reason, it is important to know which collections of measures satisfy this property. The partial reconstruction procedure provided in this paper may be used to narrow down the set of all possible measures that correspond to a given halfspace depth function. That can be used to guide the selection of an appropriate tool of depth-based analysis. The problem of determining those distributions that are uniquely characterized by their halfspace depth, however, remains open.

\appendix

\section{Proof of Theorem~\ref{theorem:support}}

For part~\ref{support2}, take $x\in\Support{\mu|_A}$ and denote $\alpha=\depth{x}{\mu}$, $D_A=\Damu\cap A$ and $U_A=\Uamu\cap A$. Because $x$ comes from the support of $\mu|_A$, we know that $\mu|_A(B_x)>0$ for any open ball $B_x$ in $A$ centered at $x$. Using Lemma~\ref{main lemma for support} we conclude that $B_x$ cannot be a subset of $\Damu\setminus \Uamu$, meaning that it also cannot be a subset of $D_A \setminus U_A$. But, because $x \in D_A \setminus U_A$, necessarily $x\in\bd_A\left(D_A\setminus U_A\right) \subseteq \bd_A\left(D_A\right)\cup\bd_A\left(U_A\right)$. Now, suppose that $x \in \bd_A\left(U_A\right)$. Then there exists a sequence $\left\{x_n\right\}_{n=1}^\infty \subset U_A \subset A$ converging to $x$. We know that $\alpha_n = \depth{x_n}{\mu} > \alpha$ for each $n$. Thus, for any $n = 1,2,\dots$ we have that $x_n \in D_{\alpha_n}(\mu)$ and $x \notin D_{\alpha_n}(\mu)$, meaning that there is a point $y_n \in A$ from the set $\bd_A\left(D_{\alpha_n}(\mu) \cap A\right)$ in the line segment $L(x,x_n)$. Since $x_n \to x$, also the sequence $\left\{y_n\right\}_{n=1}^\infty \subset \bigcup_{\beta > \alpha} \bd_A\left(D_\beta(\mu) \cap A\right)$ converges to $x$, and necessarily $x \in \cl_A\left( \bigcup_{\beta > \alpha} \bd_A\left(D_\beta(\mu) \cap A\right) \right)$ as we intended to show.

To prove part~\ref{support1}, consider $x \in \R^d$ such that $\mu(\{x\})>0$ and $\alpha = \D\left(x;\mu\right)$. Choose any $y,z\in\R^d$ such that $x\in L(y,z)$. We will prove that one of the points $y$ and $z$ has depth at most $\alpha - \mu(\{x\})$, which means that $x$ must be an extreme point of $\D_{\beta}(\mu)$ for any $\beta \in (\alpha - \mu(\{x\}), \alpha]$. Let $F \in \flag(x)$ be a minimizing flag halfspace of $\mu$ at $x$ from Lemma~\ref{theorem:Pokorny}, i.e. let $\mu(F) = \alpha$. We can write $F$ in the form of the union $\{x\}\cup\left( \bigcup_{i=1}^d\relint(H_i) \right)$ as in~\eqref{flag halfspace}. Since $H_d\in\half(x)$ is a halfspace that contains $x$ on its boundary and $x$ lies in the open line segment $L(y,z)$, one of the following must hold true with $j=d$:
\begin{enumerate}[label=(C$_{\arabic*}$), ref=(C$_{\arabic*}$)]
    \item \label{first} $L(y,z)\subset\relbd(H_j)$, or
    \item \label{second} exactly one of the points $y$ and $z$ is contained in $\relint(H_j)$.
\end{enumerate}
If~\ref{first} holds true with $j = d$, then we know that together with $x\in\relbd(H_{d-1})$ and $x\in L(y,z)$ it implies again that one of~\ref{first} or~\ref{second} must be true with $j = d-1$. We continue this procedure iteratively as $j$ decreases, until we reach an index $J\in\{1,\dots,d\}$ such that $\relint(H_J)$ contains exactly one of the points $y$ and $z$. Note that such an index must exist, since $\relbd(H_2)$ is a halfline originating at $x$, so $L(y,z)\subset \relbd(H_2)$ would imply either that $y\in\relint(H_1)$ or that $z\in\relint(H_1)$. We choose $J$ to be the largest index $j \in \{1, \dots, d\}$ satisfying~\ref{second} and assume, without loss of generality, that $y\in \relint(H_J)$. Then $L(y,z)\subset\bd(H_j)$ for each $j\in\{J+1,\dots,d)$.

Recall that for a set $A \subset \R^d$ and $u\in\R^d$ we denote by $A+u=\left\{a+u\colon a \in A\right\}$ the shift of $A$ by the vector $u$. Then for each $j\in\{1,\dots, d\}$ the $j$-dimensional halfspace $H_j+(y-x)$ satisfies $y\in\relbd(H_j + (y - x))$. Since $y\in\relint(H_J)$, it must be $\relbd\left(H_J+(y-x)\right)\subset \relint(H_J)$ and therefore 
    \begin{equation}\label{eq:aux1}
    \relint\left(H_j+(y-x)\right) \subset \left(H_j+(y-x)\right)\subset\relint(H_J)\mbox{ for }j\in\{1,\dots,J\},  
    \end{equation}
because the relative boundaries of $H_J$ and $H_J+(y-x)$ are parallel. At the same time, we have 
 \begin{equation}\label{eq:aux2}
    \relint\left(H_j+(y-x)\right)=\relint(H_j) \mbox{ for each } j\in\{J+1,\dots,d\},   
 \end{equation}
since the indices $j\in\{J+1,\dots,d\}$ all satisfy $y\in\relbd(H_j)$. Consider thus a shifted flag halfspace    
    \begin{equation}\label{eq:aux3}  
    F' = F+(y-x) = \{y\}\cup \bigcup_{j=1}^{d}\relint\left(H_j+(y-x)\right).
    \end{equation}
Using~\eqref{eq:aux1},~\eqref{eq:aux2}, and~\eqref{eq:aux3} we obtain
    \begin{equation}\label{eq:flags}
    F'\subset\{y\} \cup \relint\left(H_J\right) \cup \left( \bigcup_{j=J+1}^{d}\relint\left(H_j\right) \right) \subset F\setminus \{x\}.    
    \end{equation}
Therefore, we have $\mu(F')\leq \mu(F)-\mu(\{x\}) = \alpha - \mu(\{x\}) < \beta$ and necessarily also $\D(y;\mu) < \beta$ by Lemma~\ref{theorem:Pokorny}. Hence $y\notin \D_\beta(\mu)$. We conclude that $x \in \Damu$ cannot be contained in the relative interior of any line segment whose endpoints are both in $\D_\beta(\mu)$ for $\beta \in (\alpha - \mu(\{x\}), \alpha]$, and $x$ is therefore an extreme point of each such $\D_\beta(\mu)$. 

Consider now part~\ref{jump} and take $F$ to be any minimizing flag halfspace of $\mu$ at $x$. Then $\mu\left(F\right)=\depth{x}{\mu} = \alpha < \depth{z}{\mu}$, and necessarily $z\notin F$. Since $x\in L(y,z)$, we can use the same argumentation as in part~\ref{support1} of this proof to conclude that exactly one of the points $y$ and $z$ is contained in the relative interior of one of the closed $i$-dimensional halfspaces $H_i$ taking part in the decomposition of $F$ in~\eqref{flag halfspace}, meaning that $F\setminus\{x\}$ contains exactly one of these two endpoints $y$ and $z$. Since we found that $z\notin F$, it must be that $y\in F\setminus\{x\}$. Then from~\eqref{eq:flags} it follows that $F+(y-x)\subset F \setminus\{x\}$. Therefore,
\begin{equation}\label{eq:aux}
\mu\left(F+(y-x)\right)\leq \mu\left(F\right)-\mu\left(\{x\}\right) . \end{equation}
At the same time, Lemma~\ref{theorem:Pokorny} gives us $\D(y;\mu)\leq \mu\left(F+(y-x)\right)$, which together with~\eqref{eq:aux} finally implies
	\[	\depth{y}{\mu}\leq \mu\left(F+(y-x)\right)\leq\mu\left(F\right)-\mu\left(\{x\}\right)=\depth{x}{\mu}-\mu\left(\{x\}\right),	\] 
where the last equality follows from the fact that $F$ is a minimizing flag halfspace of $\mu$ at $x$. We proved all three parts of our main theorem.

\subsection*{Acknowledgments.}
P.~Laketa was supported by the OP RDE project ``International mobility of research, technical and administrative staff at the Charles University", grant number CZ.02.2.69/0.0/0.0/18\_053/0016976. The work of S.~Nagy was supported by Czech Science Foundation (EXPRO project n. 19-28231X).


\def\cprime{$'$} \def\polhk#1{\setbox0=\hbox{#1}{\ooalign{\hidewidth
  \lower1.5ex\hbox{`}\hidewidth\crcr\unhbox0}}}

\end{document}